\documentclass[a4paper,12pt]{article}
\usepackage[english]{babel}
\usepackage{amssymb,amsmath}

\newcommand{\mC}{\mathbb{C}}

\newcommand{\mN}{\mathbb{N}}

\newcommand{\mR}{\mathbb{R}}
\newcommand{\mT}{\mathbb{T}}
\newcommand{\mS}{\mathbb{S}}
\newcommand{\mZ}{\mathbb{Z}}

\newcommand{\te}{\theta}

\newcommand{\fhi}{\varphi}

\newcommand{\pet}{\varepsilon}

\newcommand{\cH}{\mathcal{H}}

\newcommand{\cG}{\mathcal{G}}

\newcommand{\bk}{{\bf{k}}}
\newcommand{\QED}{$\Box$}

\renewcommand{\Im}{\text{Im }}

\newcommand{\rot}{\mathrm{rot}}

\newcommand{\fm}{\phantom{-}}

\newcommand{\rmcite}[1]{{\rm{\cite{#1}}}}
\newcommand{\Spect}[1]{\mathrm{Spec}\left(#1\right)}
\newcommand{\IDS}{{\sc ids~}}
\newcommand\ii{\mathrm{i}}

\newtheorem{theorem}{Theorem}

\newtheorem{proposition}[theorem]{Proposition}
\newtheorem{lemma}[theorem]{Lemma}
\newtheorem{corollary}[theorem]{Corollary}

\newtheorem{remark}[theorem]{Remark}
\newtheorem{remarks}[theorem]{Remarks}

\title{A Nonperturbative Eliasson's Reducibility Theorem}
\author{Joaquim Puig}
\date{}
\begin{document}
\maketitle
\begin{center}
\footnotesize{Departament de Matem\`atica Aplicada I, Universitat Polit\`ecnica de Catalunya\\
Av. Diagonal 647, 08028 Barcelona, Spain \\
{\tt {joaquim.puig@upc.edu}}}
\end{center}

\abstract{This paper is concerned with discrete, one-dimensional Schr\"odinger
operators with real analytic potentials and one Diophantine frequency.
Using localization and duality we show that almost every point in the
spectrum  admits a quasi-periodic Bloch wave if the
potential is smaller than a certain constant which does not depend on the
precise Diophantine conditions. The associated
first-order system, a quasi-periodic skew-product,
is shown to be reducible for almost all values of the
energy. This is a partial nonperturbative generalization of a
reducibility theorem by Eliasson. We also 
extend nonperturbatively the genericity of Cantor spectrum
for these Schr\"odinger operators. Finally we prove that in
our setting, Cantor spectrum implies the existence of a
$G_\delta$-set of energies whose Schr\"odinger cocycle 
is not reducible to constant coefficients.


}

\tableofcontents

\section{Introduction. Main results}

Recently there has been substantial advance in the theory of
quasi-periodic Schr\"o\-din\-ger operators, both continuous and discrete,
combining spectral and dynamical techniques. These operators arise naturally in many areas of
physics and mathematics. They appear in the study of electronic properties
of solids \rmcite{aubry-andre, quasicrystals, osadchy-avron}, 
 in the theory of KdV and related equations
\rmcite{johnson-moser,johnson:sato-segal,chulaevsky}  or in Hamiltonian mechanics
\rmcite{sinai:structure}.  Moreover, their eigenvalue equations are second order 
differential or difference linear equations with quasi-periodic
coefficients like Hill's equation with quasi-periodic forcing
\rmcite{broer-puig-simo} or the Harper equation \rmcite{ketoja-satija, keller-etal} which display a
rich variety of dynamics ranging from quasi-periodicity to uniform and
nonuniform hyperbolicity. 

In this paper we pursue this fruitful combination of spectral and
dynamical methods to study discrete, one-dimensional Schr\"odinger 
operators $H_{V,\omega,\phi},$ 
\begin{equation}\label{eq:schop}
(H_{V,\omega,\phi} x)_n = x_{n+1}+x_{n-1} +  V(2\pi \omega n +
\phi) x_n, \qquad n \in \mZ,
\end{equation}
where $V:\mT \to \mR$ is a real analytic function (the \emph{potential}), $\phi
\in \mT=\mR/(2\pi\mZ)$ is a \emph{phase} and $\omega$ a
\emph{Diophantine frequency}. The latter means 
that there exist  positive constants $c$ and $\tau>1$ such that the bound
\begin{equation}\label{eq:diophantine}
\left| \sin{2\pi k  \omega}\right| > \frac{c}{|k|^\tau} 
\end{equation}
holds for any integer $k \ne 0.$ This condition will be written as $\omega
\in DC(c,\tau).$ In particular, this means that $\omega$ is \emph{nonresonant},that is
\[
\left| \sin{2\pi k  \omega}\right| \ne 0
\]
unless $k= 0.$

The best-studied example of quasi-periodic Schr\"odinger operator is the
\emph{Almost Mathieu operator} where $V(\theta) = b\cos
\theta$, being $b$ a real \emph{coupling}
parameter. A vast amount of literature is devoted to the study of the
spectral properties of this operator (see Simon
\rmcite{simon:review,simon:sXX}, Jitomirskaya \rmcite{jitomirskaya:almost,jitomirskaya:icm02}
and Last \rmcite{last:almost} for surveys and references). In \rmcite{puig,puig:thesis} 
(see also Avila \& Jitomirskaya \rmcite{avila-jitomirskaya}
for the recent extension to the remaining frequencies) 
the Cantor structure of the spectrum of the Almost Mathieu operator 
was derived from the use of 
a localization result by Jitormirskaya \rmcite{jito:metal} and a 
dynamical analysis of its eigenvalue equation, the so-called
\emph{Harper equation}. This was a long-standing conjecture
known as the ``Ten Martini Problem''. The combined approach
there is not limited to the Almost Mathieu operator as we 
plan to make evident in this paper. 

The eigenvalue equation of a quasi-periodic Schr\"odinger operator
$H_{V,\omega,\phi}$ is the following Harper-like equation
\begin{equation}\label{eq:hill} 
x_{n+1}+x_{n-1} + V(2\pi \omega n + \phi) x_n= a x_n, \qquad n
\in \mZ,
\end{equation}
where $a\in \mR$ is called the  \emph{energy} or \emph{spectral parameter}. Since
we want to study dynamical properties of this equation, it is better to
transform it into a first-order system, obtaining the associated \emph{quasi-periodic skew-product} on $\mR^2 \times \mT$,
\begin{equation}
\left(\begin{array}{c}
x_{n+1} \\
x_{n} 
\end{array}\right)=
\left(
\begin{array}{cc}
a - V(\te_n)  & -1 \\
1 & \phantom{-}0
\end{array}\right)
\left(\begin{array}{c}
x_{n} \\
x_{n-1} 
\end{array}\right),\qquad \te_{n+1} = \te_n + 2\pi \omega,
\end{equation}
which can be seen as an iteration of the corresponding \emph{Schr\"odinger
cocycle}, $(A_{a,V},\omega),$ on $SL(2,\mR) \times \mT$.
Here $A_{a,V}$ denotes the matrix-valued function
\begin{equation}\label{eq:schrodinger_cocycle}
A_{a,V}(\theta)= \left(
\begin{array}{cc}
a -  V(\theta)  & -1 \\
1 &  \phantom{-}0
\end{array}
\right), \qquad \theta \in \mT.
\end{equation}

The simplest class of Schr\"odinger cocycles occurs when
$V=0$ because in this case $A_{a,V}$ does not depend on
$\theta$ (we will say that the corresponding cocycle is in
\emph{constant coefficients}). In analogy with periodic
differential equations we may try to reduce a
quasi-periodic Schr\"odinger cocycle to constant
coefficients.  Let us introduce first the notion of
conjugation between cocycles, not necessarily of
Schr\"odinger type. Here we restrict ourselves to the case of $SL(2,\mR)$-valued cocycles, although the notion of reducibility applies to more general cocycles (see  \rmcite{puig:thesis} for an exposition). Two cocycles 
$(A,\omega)$ and $(B,\omega)$ on $SL(2,\mR)\times \mT$ 
 are \emph{conjugated}
if there exists a continuous and nonsingular \emph{conjugation}
$Z:\mT \to SL(2,\mR)$ such that the relation
\[
A(\theta) Z(\theta) = Z(\theta + 2\pi \omega) B(\theta), \qquad \theta
\in \mT.
\]
holds for all $\theta \in \mT.$ In this case the corresponding 
quasi-periodic skew-products
\[
u_{n+1} =  A(\theta_n) u_n,\qquad \theta_{n+1 }= \theta_n + 2\pi \omega
\]
and
\[
v_{n+1} =  B(\theta_n) v_n,\qquad \theta_{n+1 }= \theta_n + 2\pi \omega
\]
are conjugated through the change of variables $u=Zv$, so that they share
the same dynamical properties.

A cocycle $(A,\omega)$ is
\emph{reducible to constant coefficients} if it is conjugated to a cocycle
$(B,\omega)$ with $B$ not depending on $\theta$ (i.e. with
constant coefficients). In this case, $B$ is
called a \emph{Floquet} matrix. Sometimes it may be necessary to ``halve
the frequency'' if we do not want to complexify the system
(although this case will not be treated in this paper).
In contrast with the situation in the periodic case 
(when the nonresonance condition fails), quasi-periodic cocycles need
not to be reducible to constant coefficients (see Theorem
\ref{res:nored} and the following Remark
\ref{rem:nored}).

The reducibility of a Schr\"odinger cocycle and the eigenvalues of the
reduced Floquet matrix have implications for the spectrum of the corresponding
Schr\"odinger operator. The spectrum of
$H_{V,\omega,\phi}$ on $l^2(\mZ)$ is a compact subset of
the real line which we denote by $\sigma(V,\omega)$ since
it does not depend on $\phi$. It is known that an
energy $a$ lies in the spectrum of a Schr\"odinger operator if, and only
if, the corresponding skew-product has an \emph{exponential dichotomy} (it is
\emph{uniformly hyperbolic}), see Johnson \rmcite{johnson:recurrent}. Under our assumptions, $V$ real analytic and
$\omega$ Diophantine, this is equivalent to the reducibility of the
Schr\"odinger cocycle to constant coefficients with a hyperbolic Floquet
matrix (all its eigenvalues are outside the unit circle), see 
Johnson \rmcite{johnson:analyticity}.

For energies in the spectrum, the situation is much more
involved. However, when in addition to the present hypothesis, the 
potential $V$ is small (in some complex neighbourhood
around $\mT$ which depends on $c$ and $\tau$), then reducibility
can be obtained by {\sc{kam}} methods for a set of energies
in the spectrum of large measure, see Dinaburg \& Sinai
\rmcite{dinaburg-sinai} and Moser \& P\"oschel \rmcite{moser-poschel}.

A breakthrough in the {\sc{kam}} approach came with Eliasson 
\rmcite{eliasson:floquet} who proved, 
among other statements, that reducibility to constant coefficients 
holds for almost every energy provided the potential is
small enough and $\omega$ is Diophantine,  $\omega \in
DC(c,\tau)$ for some $c$ and $\tau$ (see Section
\ref{sec:eliasson} for a more precise formulation.
Like the results in the previous paragraph, the smallness
condition here depends on the precise 
Diophantine conditions on $\omega$. Eliasson's result, as well as the above
{\sc{kam}} results presented above, holds for real analytic
or $C^\infty$ potentials $V:\mT \to \mR$ depending on several frequencies. 

Our main result states that in the presence of only one 
frequency, $d=1$ in the notation above,  
the smallness condition in Eliasson's theorem does not depend on the 
constants $c$ and $\tau$ of the Diophantine condition as
long as $\omega$ is Diophantine (that is, it is ``nonperturbative'' in some sense).
To be more precise, we consider real analytic potentials $V:\mT \to \mR$ 
having an analytic extension to $|\Im \theta|<\rho,$ for some $\rho>0$ such
that
\[
|V|_{\rho} := \sup_{|\Im \theta|<\rho} |V(\theta)| <\infty,
\]
(the set of such potentials will be denoted by $C^a_\rho(\mT,\mR)$) and a
Diophantine frequency $\omega \in DC(c,\tau)$ for some positive
constants $c$ and $\tau>1$. Our extension of Eliasson's theorem reads as
follows.

\begin{theorem}\label{res:main}
Let $\rho>0$ be a positive number. Then, there is a constant
$\pet_0=\pet_0(\rho)$ such that, for any real analytic $V \in
C^a_\rho(\mT,\mR)$ with
\[
\left|V \right|_{\rho} < \pet_0,
\]
the Schr\"odinger cocycle $(A_{a,V},\omega)$ is reducible
to constant coefficients for every Diophantine frequency
$\omega$ and almost all $a \in \mR$ (with respect to Lebesgue measure).
\end{theorem}
The proof of this Theorem will be given in Section
\ref{sec:proofs}.

\begin{remarks}
\mbox{}
\begin{enumerate}
\item Recently  Avila \& Krikorian \rmcite{avila-krikorian}
      proved  Theorem \ref{res:main} with more
      restrictive hypothesis on $\omega$ (although it is also a full
      measure condition). In fact, we will see
      that both results follow from a nonperturbative theorem
      on localization by Bourgain \& Jitomirskaya
      \rmcite{bourgain-jitomirskaya:absolutely}.
\item When the potential is defined on a $d$-dimensional torus, $V: \mT^d
      \to \mR$, Eliasson's theorem holds, but a nonperturbative version
      like Theorem \ref{res:main} cannot be true, as Bourgain showed in
      \rmcite{bourgain:spectrum2}.
      Indeed, he proved
      that,  if $V:\mT^2\to  \mR$ is a trigonometric polynomial with a nondegenerate
      maximum, there is a set of frequencies $\omega \in \mR^2$, with
      positive Lebesgue measure, for which the operators
      $H_{V,\omega,\phi}$ have some point spectrum. This point
      spectrum is incompatible with reducibility to constant
      coefficients. See  Bourgain
      \rmcite{bourgain:quasi-periodic,bourgain:book} for the
      differences  between the cases of one and several frequencies.
\item Let us stress that Theorem \ref{res:main} is not
      a full nonperturbative version of Eliasson's theorem because
      the set of energies whose corresponding
      Schr\"odinger cocycle is reducible to constant coefficients
      is not explicitly characterized as it is in Eliasson's theorem (see
      Section \ref{sec:eliasson}). 
\item The smallness condition in Theorem \ref{res:main},
      $\pet_0=\pet_0(\rho)$ is given by the localization
      result in \rmcite{bourgain-jitomirskaya:absolutely}
      and can be explicitly given in terms of $\|V\|_1$,
      $\|V\|_2$, $\|V\|_\infty$ and $\rho$. If $V$ is kept
      fixed, then $\pet_0(\rho)=O(\rho)$ as $\rho\to 0.$
\end{enumerate}
\end{remarks}

Theorem \ref{res:main} has some implications for the spectral
properties of Schr\"odinger operators. The first one refers
to the properties of solutions of the eigenvalue equation.
An immediate application of Theorem \ref{res:main} is the
existence of analytic quasi-periodic \emph{Bloch waves}\index{Bloch wave} 
for almost all $a$ in the spectrum. An analytic quasi-periodic 
Bloch wave for a Harper-like equation (\ref{eq:hill}) is 
a solution of the form
\begin{equation}\label{eq:blochwave}
x_n(\phi) = e^{\ii \fhi n} f\left(2\pi\omega n +
\phi\right) , \qquad n \in \mZ,
\end{equation}
where $\fhi \in [0,2\pi)$ is called the \emph{Floquet exponent} and
$f: \mT \to \mC$ is a nontrivial analytic function. In Section 
\ref{sec:proofs} we will also prove the following about 
the existence of quasi-periodic Bloch waves.

\begin{corollary}\label{res:bloch}
Let $\rho$, $\pet_0$ $V$ and $\omega$  be as in Theorem
\ref{res:main}. Then for (Lebesgue) almost all values of $a$ in
the spectrum $\sigma(V,\omega)$, the equation
\[
x_{n+1}+x_{n-1} + V(2\pi \omega n + \phi) x_n= a x_n,
\]
has analytic quasi-periodic Bloch waves.
\end{corollary}

Our second application deals with the structure of the
spectrum of quasi-periodic Schr\"odinger operators. 
The Cantor structure of the spectrum is not specific of the
Almost Mathieu operator. Indeed, although the proof of the
``Ten Martini Problem'' \rmcite{puig} is restricted to this
model, Cantor spectrum is ``generic'' in our setting.
More precisely, if we consider 
the set $C^a_\rho(\mT,\mR)$ of real analytic functions furnished with the
$|\cdot|_\rho$ norm, one has the following.

\begin{theorem}\label{res:cantor}
Let $\rho>0$. Then there is a constant $\pet=\pet(\rho)$ such that for
every Diophantine $\omega$ there is a generic set of real analytic potentials $V \in
C^a_\rho(\mT,\mR)$ with  $|V|_\rho<\pet$ such that the
spectrum of  the Schr\"odinger operator $H_{V,\omega,\phi}$
is a Cantor set.
\end{theorem}

Our last application is concerned with the existence of
Schr\"odinger cocycles which are not reducible to constant
coefficients. We will see that Cantor spectrum and
nonreducibility are related concepts in our setting.

\begin{theorem}\label{res:nored}
Let $\rho>0$. Then there is a constant $\pet=\pet(\rho)$
such that if $\omega$ is  Diophantine and $V \in
C^a_\rho(\mT,\mR)$, with  $|V|_\rho<\pet$, is such that 
$\sigma(V,\omega)$ is a Cantor set then for a $G_\delta$-dense 
subset of energies in the spectrum the
corresponding Schr\"odinger cocycle is 
not reducible to constant coefficients (by a continuous
transformation).
\end{theorem}

\begin{remark}
According to Theorem \ref{res:main}, the $G_\delta$-set above has zero Lebesgue measure.
\end{remark}

In particular, using Theorem \ref{res:cantor} it is
possible to give a nonperturbative version of a result in
Eliasson \rmcite{eliasson:floquet}, namely, that the
existence of a $G_\delta$-subset of ``nonreducible
energies'' is  a generic property. 

\begin{remark}\label{rem:nored}
This kind of nonreducibility holds for a zero measure subset
of energies and the corresponding Schr\"odinger cocycle has
zero Lyapunov exponent. Nonreducibility results can be
obtained using, for instance, a result by Sorets \& Spencer
\rmcite{sorets-spencer}, who prove that if
the potential $V$ is large enough then the Lyapunov exponent
is positive for all energies in the spectrum and this
prevents reducibility.
\end{remark}

Let us finally outline the contents of this paper. In Section
\ref{sec:prelim} we introduce some of the preliminaries
needed for the proof of the main theorem. In 
Section \ref{sec:proofs} this is
used to prove \ref{res:main} and \ref{res:bloch} using a
similar technique to the one used for the Almost Mathieu
operator. The applications are included in Section
\ref{sec:applications}.

\section{Preliminaries}\label{sec:prelim}

In this section we present some of the tools that will be
needed in the proof of Theorem \ref{res:main}. As said in the
introduction, we plan to extend some of the ideas in the proof of
the ``Ten Martini Problem'' given in \rmcite{puig} with the aid of a 
nonperturbative localization result for long-range potentials by 
Bourgain \& Jitomirskaya \rmcite{bourgain-jitomirskaya}. In
Section \ref{sec:aubry} we introduce a convenient version of Aubry duality, 
which will lead us to consider certain long-range operators which 
are not of Schr\"odinger type. In Section \ref{sec:ids} we
give the definition and some properties of the integrated
density of states \IDS for these operators and its relation with Aubry duality.
Finally, in Section \ref{sec:eliasson} this \IDS is
linked to the fibered rotation number of quasi-periodic
cocycles to give a more precise version of Eliasson's
result.

\subsection{Aubry Duality}\label{sec:aubry}

Aubry Duality \rmcite{aubry-andre} was originally introduced for the study 
of the Almost Mathieu operator but the same idea (which is
Fourier transform) works for other potentials. Let us give
first the heuristic approach and then a more rigorous one.

Assume that 
\[
x_{n+1}+x_{n-1} + V(2\pi \omega n + \phi) x_n= a x_n,
\]
has an analytic quasi-periodic Bloch wave,
\begin{equation}\label{eq:bloch_wave}
x_n = e^{\ii \fhi n} {\tilde \psi}\left(2\pi\omega
n + \phi \right),
\end{equation}
being ${\tilde \psi} : \mT \to \mC$ analytic and $\fhi
\in [0,2\pi)$ the Floquet exponent. If $(\psi_n)_{n\in
\mZ}$ are the Fourier
coefficients of ${\tilde \psi}$, a computation shows that
they satisfy the following difference equation
\[
\sum_{k \in \mZ} V_k \psi_{n-k} + 2 \cos\left(2\pi\omega
n + \fhi \right)
\psi_n = a \psi_n \qquad n \in \mZ,
\]
where $(V_k)_{k\in\mZ}$ are the Fourier coefficients of $V$,
\[
V(\theta)= \sum_{k\in\mZ} V_k e^{\ii k \theta}.
\]
This difference equation is the eigenvalue equation of the operator $L_{V,\omega,\fhi}$
\[
\left(L_{V,\omega,\fhi} \psi \right)_n = 
 \sum_{k \in \mZ} V_k \psi_{n-k} + 
2\cos\left(2\pi\omega n + \fhi \right)
\psi_n
\]
which we call  a \emph{dual operator} of $H_{V,\omega,\phi}$.
This is a self-adjoint and bounded operator on $l^2(\mZ)$ (because
$V$ is real analytic)  but  it is
not a Schr\"odinger operator unless $V$ is exactly the cosine (this is what 
makes the Almost Mathieu operator so special).
Such an operator will be called a \emph{long-range
(quasi-periodic) operator} even if it may be a
finite-differences operator (if $V$ is a trigonometric
polynomial). 

If $\omega$ is nonresonant, the spectrum of the long-range
operators $L_{V,\omega,\fhi}$ does not depend on the chosen
$\fhi$, so that one can write
\[
\sigma^L(V,\omega)= \Spect{L_{V,\omega,\fhi}}.
\]

This naive approach to Aubry duality shows that whenever $a$ is a
value in the spectrum $\sigma^H(V,\omega)$ such that $(x_n)_{n\in
\mZ}$ is an analytic quasi-periodic Bloch wave 
with Floquet exponent $\fhi$, then $a$ is a
point eigenvalue of the dual operator $L_{V,\omega,\fhi}$ whose
eigenvector decays exponentially and, thus,
$a\in\sigma^L(V,\omega)$. The converse is also true: one
can pass from exponentially decaying eigenvalues of
$L_{V,\omega,\fhi}$ to quasi-periodic Bloch waves of
$H_{V,\omega,\phi}$ with Floquet exponent $\fhi$.

The argument given above heavily relies on the existence of
quasi-periodic Bloch waves or, equivalently, exponentially
localized eigenvectors. Nevertheless both operators
can be related without the assumption of such point
eigenvalues. This was done by Avron \& Simon
\rmcite{avron-simon}. Here we will follow the idea by Gordon,
Jitomirskaya, Last \& Simon
\rmcite{gordon-jitomirskaya-last-simon} (see also Chulaevsky \& Delyon \rmcite{chulaevsky-delyon}),
who studied duality for the Almost Mathieu operator, although it can be extended to
the general case, see  Bourgain \& Jitomirskaya \rmcite{bourgain-jitomirskaya}. The
idea is to shift to  more general spaces where the extensions of
quasi-periodic Schr\"odinger operators and their duals are unitarily equivalent.
Note that it is not true that the operators $H_{V,\omega,\phi}$
and $L_{V,\omega,\fhi}$ are unitarily equivalent, since their
spectral measures will, in general, be very different.

Let us consider the following Hilbert space,
\[
\cH = L^2\left(\mT \times \mZ \right),
\]
which consists of functions $\Psi=\Psi(\theta,n)$ satisfying
\[
\sum_{n\in \mZ} \int_\mT \left| \Psi(\theta,n) \right|^2 d\theta
< \infty.
\]

The extensions of the Schr\"odinger operators $H$ and their
long-range duals $L$ to $\cH$ are given in terms of their
\emph{direct integrals}\index{direct integral of a quasi-periodic operator},
which we now define. The  \emph{direct integral} of the Schr\"odinger operator
$H_{V,\omega,\phi}$, is the operator ${\tilde H}_{V,\omega}$, defined as
\[
\left({\tilde H}_{V,\omega} \Psi \right) (\theta,n) =
\Psi(\theta,n+1) + \Psi(\theta,n-1) + V(2\pi\omega n+ \theta)
\Psi(\theta,n),
\]
and the direct integral of $L_{V,\omega,\fhi}$, denoted as ${\tilde
L}_{V,\omega}$, is
\[
\left({\tilde L}_{V,\omega} \Psi \right) (\theta,n) =
 \sum_{k \in \mZ} V_k \Psi (\theta,n-k) + 2
\cos\left(2\pi\omega n + \theta \right)
\Psi(\theta,n).
\]
These two operators are bounded and self-adjoint in $\cH$.
Let us now see that, for any fixed real analytic $V$ and
nonresonant frequency $\omega$, the direct integrals ${\tilde
H}_{V,\omega}$ and ${\tilde L}_{V,\omega}$ are unitarily
equivalent; i.e. there exists a unitary operator $U$
on $\cH$ such that the conjugation
\[
{\tilde H}_{V,\omega} U = U {\tilde L}_{V,\omega} 
\]
holds. By analogy with the heuristic approach to Aubry duality
in the beginning of this section, let $U$ be the following
operator on $\cH$,
\[
\left(U \Psi \right)(\theta,n)= {\hat \Psi} \left(n, \theta 
+ 2\pi\omega n  \right),
\]
where $\hat \Psi$ is the Fourier transform. At a formal level this 
acts as
\[
\left(U \Psi\right)(\theta, n)= 
\sum_{k \in\mZ} \int_{\mT}
\Psi(\phi,k) e^{-\ii n\phi} e^{-\ii\left(\theta +2\pi\omega n \right)k} d\theta
\]
if we disregard the convergence of the sum in $k.$ The map $U$ is unitary and satisfies
\[
{\tilde H}_{V,\omega} U = U {\tilde L}_{V,\omega}
\]
by construction of the dual long-range operators in terms of the
Schr\"o\-din\-ger operators. Therefore, the direct integrals ${\tilde
H}_{V,\omega}$ and ${\tilde L}_{V,\omega}$ are unitarily
equivalent and, in particular, their spectra are the same,
\begin{multline}
\sigma^H(V,\omega)= 
\bigcup_{\phi \in \mT} \Spect{H_{V,\omega,\phi}}
= \Spect{{\tilde H}_{V,\omega}}= \\
\Spect{{\tilde L}_{V,\omega}}=
\bigcup_{\fhi \in \mT} \Spect{H_{V,\omega,\fhi}}= 
\sigma^L(V,\omega).\nonumber
\end{multline}
Hence, the spectrum of a quasi-periodic Schr\"odinger operator and
its dual are the same. In the next section we will introduce the
integrated density of states for Schr\"odinger operators
(and their long-range duals) and we will see
that this function is  preserved by Aubry duality.

\subsection{The integrated density of states and duality}\label{sec:ids}

The integrated density of states, \IDS for short, is a very useful
object for the description of the spectrum of quasi-periodic Schr\"odinger
ope\-ra\-tors and more general quasi-periodic self-adjoint operators. 
Here we want to introduce it both for quasi-periodic Schr\"odinger operators
and their long-range duals. In order to give a unified
approach, let us consider a more general class of
operators. 

If $V,W:\mT \to \mR$ are
real analytic functions, $(V_k)_k$ and $(W_k)_k$ are their
Fourier coefficients, $\omega$ is a nonresonant frequency 
and $\phi \in \mT$, let
$K_{W,V,\omega,\phi}$ be the following operator
\[
\left(K_{W,V,\omega,\phi} x\right)_n =
\sum_{k \in \mZ} W_k x_{n-k} + V(2\pi\omega n + \phi) x_n
\]
acting on $l^2(\mZ),$ which is bounded and self-adjoint.
The operators in the previous section occur as particular cases,
\[
H_{V,\omega,\phi} = K_{2\cos ,V,\omega,\phi} \quad\text{ and }\quad
L_{V,\omega,\phi} = K_{V,2\cos,\omega,\phi}.
\]

Let us now define the \IDS for the operators $K_{W,V,\omega,\phi}.$
Take some integer $N>0$ and consider $K_{W,V,\omega,\phi}^N$, the
restriction of the operator $K_{W,V,\omega,\phi}$ to the interval $[-N,N]$ with zero boundary conditions. Let 
\[
k^N_{W,V,\omega,\phi}(a) = \frac{1}{2N+1} \# \left\{\text{eigenvalues}
\le a \text{ of } K_{W,V,\omega,\phi}^N \right\}.
\]
Then, due to the nonresonant character of $\omega,$ the limit
\[
\lim_{N\to\infty} k^N_{W,V,\omega,\phi}(a) 
\]
exists, it is independent of $\phi$ and of the boundary
conditions imposed above. It is called the \emph{integrated
density of states} of the operator $K_{W,V,\omega,\phi}.$
We will write this as $k_{a,W,V,\omega}(a).$  The map
\begin{equation}\label{eq:idsmap}
a \in \mR \mapsto k_{a,W,V,\omega}(a)
\end{equation}
is increasing and it is constant exactly at the open
intervals in the resolvent set of the spectrum of 
$K_{W,V,\omega,\phi}$. It is the distribution function of a
Borel measure $n_{W,V,\omega}$, 
\[
k_{a,W,V,\omega}(a)= \int_{-\infty}^{a} dn_{W,V,\omega}(\lambda)
\]
called is the \emph{density of states} of the operator
$K_{W,V,\omega,\phi},$
which is supported on the spectrum of  $K_{W,V,\omega,\phi}$.
In the Schr\"odinger case we will use the notations
\[
k^H_{V,\omega}(a)=k_{2\cos,V,\omega}(a),\quad n^H_{V,\omega}=
n_{2\cos,V,\omega}
\]
and 
\[
k^L_{V,\omega}(a)=k_{V,2\cos,\omega}(a),\quad n^L_{V,\omega}=
n_{V,2\cos,\omega}.
\]
for their long-range duals.

The \IDS of the operators $K_{W,V,\omega,\phi}$ can be seen as an average in
$\phi$ of the spectral  measures of the operators (see  Avron \& Simon \rmcite{avron-simon}).
By the spectral theorem we know that there is a Borel measure $\mu_{\phi}$ such
that
\begin{equation}\label{eq:papers-*}
\langle \delta_0, f\left(K_{W,V,\omega,\phi}\right) \delta_0 \rangle_{l^2(\mZ)} = 
\int f(\lambda) d\mu_\phi(\lambda)
\end{equation}
for every continuous function $f$, being $\delta_0$ the delta function.
The measures $\mu_{\phi}$ are spectral measures in the sense that the spectral
projection of $K_{W,V,\omega,\phi}$ over a certain subset $A$ of the spectrum is
zero if, and only if, $\mu_{\phi}(A)=0$. Avron \& Simon prove that, for any
continuous function $f$
\[
\int f(\lambda) dn^K_{W,V,\omega}(\lambda) 
= \int_{\mT} d\phi \int f(\lambda) d\mu_{\phi}.
\]

An approximation argument shows that, for any Borel subset of the spectrum, $A
\subset \sigma(K_{W,V,\omega}),$
\[
n^K_{W,V,\omega}(A)= \int_{\mT} \mu_{\phi}(A) d\phi.
\]
In particular, $n^K_{W,V,\omega}(A)=0$ if $\mu_{\phi}(A)=0$ for Lebesgue almost
every $\phi \in \mT$. Using this characterization of the \IDS one can
prove the following adaption of the duality of the \IDS given in
\rmcite{gordon-jitomirskaya-last-simon}.

\begin{theorem}[\rmcite{gordon-jitomirskaya-last-simon}]
Let $k^L_{V,\omega}$ and $k^H_{V,\omega}$ be the integrated density of
states of $H_{V,\omega,\phi}$ and $L_{V,\omega,\fhi}$
respectively, for some real analytic $V:\mT \to \mR$ and
nonresonant frequency $\omega$. Then
\[
k^L_{V,\omega}(a) = k^H_{V,\omega}(a)
\]
for all $a \in \mR$. 
\end{theorem}

\noindent{\bf{Proof:}} Let  
\[
g(\theta,n)= \delta_{n,0}
\]
which belongs to $\cH$. Then $U g= g.$ Moreover by
(\ref{eq:papers-*}) and the
unitary equivalence between $\tilde{H}_{V,\omega}$ and $\tilde{L}_{V,\omega}$
we have that, for any continuous $f,$
\begin{multline}
\langle g, f\left({\tilde H}_{V,\omega}\right) g\rangle_{\cH} =
\langle U g, U f\left({\tilde H}_{V,\omega}\right)
g\rangle_{\cH} = \\
\langle U g, U f\left({\tilde H}_{V,\omega}\right) U^{-1} U g\rangle_{\cH} =
\langle g, f\left({\tilde L}_{V,\omega}\right)
g\rangle_{\cH}. \nonumber
\end{multline}
Therefore, since $n^{H}_{V,\omega}$ and $n^{L}_{V,\omega}$ are the Borel
measures such that
\[
\langle g, f\left({\tilde H}_{V,\omega}\right) g\rangle_{\cH} = \int f(\lambda)
d n^{H}_{V,\omega}(\lambda)
\]
and 
\[
\langle g, f\left({\tilde L}_{V,\omega}\right) g\rangle_{\cH} = \int f(\lambda)
d n^{L}_{V,\omega}(\lambda)
\]
for every continuous $f$ the two measures must coincide
(and also their distribution functions, $k^L_{V,\omega}$
and $k^H_{V,\omega}$).\hfill \QED

Let us end this section summing up some facts 
useful in the sequel.

\begin{proposition}\label{res:zeromeasure}
Let $V$ be real analytic, $\omega$ nonresonant and
$\mu_\phi$ a spectral measure of $L_{V,\omega,\phi}$. Assume
that there is a measurable set $A$ such that 
\[
\mu_\phi(A)= 0
\]
for almost every $\phi \in \mT$. Then $n_{V,\omega}^L(A)=0$ and
$n_{V,\omega}^H(A)=0.$ 
\end{proposition}

\subsection{The rotation number and Eliasson's theorem
revisited}\label{sec:eliasson}

We have seen in the previous section that it is
possible to assign an \IDS for quasi-periodic
Schr\"odinger cocycles using its associated
operator. Here we will see
that it is possible to define an extension of this
object, the fibered rotation number, for more
general quasi-periodic cocycles. This object,
introduced originally by Herman \rmcite{herman} in
this discrete case (see also Johnson \& Moser
\rmcite{johnson-moser}, Delyon \& Souillard
\rmcite{delyon-souillard}), allows us to give a
version of Eliasson's theorem for these cocycles. Let us
follow the presentation by Krikorian \rmcite{krikorian:sl2}.

Let $(A,\omega)$ be a quasi-periodic cocycle on
$SL(2,\mR)\times\mT$ which is \emph{homotopic to the
identity}. That is  $A:\mT \to SL(2,\mR)$ is a continuous
map (although we will later assume that it is real
analytic) that is homotopic to the identity (for 
example a Schr\"odinger cocycle). The fibered
rotation number, which we now introduce measures how
solutions wind around the origin in $\mR^2$ in average.

Let $\mS^1$ be the set of unit vectors of $\mR^2$ and let us
denote by $p:\mR \to \mS^1$ the projection given by the
exponential $p(t)= e^{\ii t}$, identifying $\mR^2$ with $\mC$.
Because of the linear character of the cocycle and 
the fact that it is homotopic to the identity, the continuous map
\begin{equation}
\begin{array}{rrcl}
F: & \mS^1 \times \mT  & \longrightarrow     & \mS^1 \times \mT \\
      & (v,\theta)    & \mapsto & \left(
      \frac{A(\theta)v}{\| A(\theta) v \|}, \theta+2\pi\omega \right)
\end{array}
\end{equation}
is also homotopic to the identity. Therefore, it
admits a continuous lift
$\tilde{F}: \mR \times \mT \to  \mR \times \mT$ of the form:
\[
\tilde{F}(t,\theta)= \left( t+ f(\theta,t), \theta + 2\pi\omega\right)
\]
such that
\[
f(t +2\pi,\theta + 2\pi)= f(t,\theta) \text{ and } p\left(t + f(t,\theta)\right)=
\frac{A(\theta)p(t)}{\| A(\theta) p(t) \|}
\]
for all $t \in \mR$ and $\theta \in \mT$. The map $f$ is independent of
the choice of $\tilde{F}$ up to the addition of a
constant $2\pi k$, with $k \in \mZ$. Since the iteration
$\theta \mapsto \theta + 2\pi\omega$ is
uniquely ergodic on $\mT$ for all $(t,\theta) \in
\mR\times\mT$, one has that  the limit
\[
\lim_{N \to \infty} \frac{1}{2\pi N} \sum_{n=0}^{N-1}
f\left(\tilde{F}^n(t,\theta)\right)
\]
exists modulus $\mZ$ and it is independent of
$(t,\theta)$, see Herman
\rmcite{herman}. This object is called the
\emph{fibered rotation number} of $(A,\omega)$, 
and it will be denoted by
$\rot_f(A,\omega)$. The fibered rotation number of a Harper-like
equation is defined as the fibered rotation number of the 
associated Schr\"odinger cocycle on $SL(2,\mR) \times \mT$ 
and will be denoted as $\rot_f(a,V,\omega)$.

The rotation number of a Harper-like equation can be linked
to its \IDS. Indeed, using a suspension argument (see Johnson
\rmcite{johnson:review}) it can be seen that
\[
\rot_f(a,V,\omega)= \frac{1}{2}k_{V,\omega}(a) \qquad
(\text{mod. } \mZ). 
\]

The rotation number is not invariant under conjugation, but
one has the following.

\begin{proposition}[cf. \rmcite{krikorian:sl2}]
Let $\omega$ be nonresonant and $(A_1,\omega)$ and
$(A_2,\omega)$ be  two quasi-periodic cocycles
on $SL(2,\mR) \times \mT$ homotopic to the identity. If they
are conjugated for some continuous $Z:\mT \to SL(2,\mR),$ then 
\[
\rot_f(A_1,\omega)= \rot_f(A_2,\omega) + 
\langle \bk, \omega \rangle  \text{ modulus } \mZ,
\]
where $\bk \in \mZ$ is the degree of the map $Z:\mT \to SL(2,\mR)$.
If the conjugation $Z$ is not defined on $\mT$ but on
$(\mR/(4\pi\mZ))$ and it has degree $\bk \in \mT$, then
\[
\rot_f(A_1,\omega)= \rot_f(A_2,\omega) +
\frac{1}{2}\langle \bk,\omega
\rangle.
\]
\end{proposition}

Keeping in mind this result, we can define two classes of 
rotation numbers which are preserved under conjugation. An
important class is that of resonant rotation numbers. A number of the form
\[
\alpha = \frac{1}{2}\langle \bk, \omega
\rangle,\qquad  (\text{mod } \frac{1}{2}\mZ) 
\]
for some $k\in \mZ$ is called \emph{resonant with respect to} $\omega$. We can also define the class of fibered rotation numbers which are
{Diophantine with respect to} $\omega$. Its elements are the numbers 
$\alpha$ such that the bound
\[
\left|\sin\left(\pi\left(2\alpha - \langle
k,\omega\rangle\right)\right)\right|
\ge \frac{K}{|k|^\sigma},
\]
holds for all $k \in \mZ-\{0\}$ and suitable fixed positive
constants $K$ and $\sigma$. Both classes of rotation
numbers are constant under conjugation.

With these definitions  we can give a more precise
version of Eliasson's reducibility theorem for general
quasi-periodic cocycles on $SL(2,\mR) \times \mT$
homotopic to the identity. Again the result is valid
for more than one frequency, but we restrict
ourselves to this one-dimensional case.

\begin{theorem}[\rmcite{eliasson:floquet}]\label{res:eliasson} 
Let $\rho>0$, $\omega \in DC(c,\tau)$ be Diophantine and $A_0$ be a
matrix in $SL(2,\mR)$. Then there is a
constant $C= C(c,\tau,\rho,|A_0|)$ such that, if 
$A \in C^a_\rho(\mT,SL(2,\mR))$ is  real analytic
with
\[
\left| A - A_0 \right|_{\rho} < C
\]
and the rotation number of the cocycle $(A,\omega)$ is either
resonant or Diophantine with respect to $\omega$,
then $(A,\omega)$ is reducible to constant
coefficients of a quasi-periodic (perhaps with frequency $\omega/2$) 
and analytic transformation. 
\end{theorem}

\begin{remark}
The proof of this theorem was  originally given in
\rmcite{eliasson:floquet} in the continuous case and for
Schr\"odinger operators (instead of cocycles), although it
extends to the setting of Theorem \ref{res:eliasson}.
\end{remark}

Applied to Schr\"odinger cocycles one obtains the
perturbative version of Theorem \ref{res:main} with
the additional characterization of the set of
reducible energies in terms of its rotation number.
More precisely, the theorem above implies that the set of
``reducible'' rotation numbers is of full measure in
$\mT$. To obtain a full-measure condition on the
energies it is necessary to use some facts on the
growth of the rotation number at these reducible
points which will be also used in Section
\ref{sec:proofs} and which are due to Deift \& Simon
\rmcite{deift-simon}.

\section{Proof of Theorem \ref{res:main}}\label{sec:proofs}

We are now ready to show  that Theorem \ref{res:main} is a
direct consequence of the following result by Bourgain \&
Jitomirskaya \rmcite{bourgain-jitomirskaya}, which we restate in a
convenient way:

\begin{theorem}[\rmcite{bourgain-jitomirskaya}]\label{res:bourgain-jitomirskaya}
Let $\rho>0$ be a positive number. Then there is a constant
$\pet_0=\pet_0(\rho)$ such that, for any real analytic $V \in
C^a_\rho(\mT,\mR)$ with 
\[
\left|V \right|_{\rho}< \pet_0,
\]
and Diophantine $\omega$ there is a set $\Phi \subset \mT$, of
zero (Lebesgue) measure such that, if $\phi\not\in\Phi$, the
operator $L_{V,\omega,\phi}$ has pure point spectrum with
exponentially decaying eigenfunctions.
\end{theorem}
\smallskip

\begin{remarks}
\mbox{}
\begin{enumerate}
\item In \rmcite{bourgain-jitomirskaya:absolutely}, the bound
   $\pet_0$ depends on $\|V\|_{1},$ $\|V\|_{2},$ $\|V\|_{\infty}$
   and $\rho$. Since $V$ belongs to $C^a_\rho(\mT,\mR)$, all
   these previous norms can be controlled by $|V|_{\rho}$.
\item  The set $\Phi$ consists of those phases $\phi$ for which the
    relation
    \begin{equation}\label{eq:resonant}
    \left| \sin \left( \phi + \pi k \omega  \right) \right| <
    \exp \left( - |k|^{\frac{1}{2\tau}}
    \right) 
    \end{equation}
    holds for infinitely many values of $k$, where
    $\omega \in DC(c,\tau)$. For any Diophantine $\omega$,
    this is a set of zero Lebesgue measure.
\end{enumerate}
\end{remarks}

Our strategy to prove Theorem \ref{res:main} will
be, first of all, to show that Corollary
\ref{res:bloch} is a simple consequence of Theorem
\ref{res:bourgain-jitomirskaya} and the duality of
the \IDS. Then, in Section \ref{sec:bloch2red} it
will be shown that Corollary \ref{res:bloch} actually implies
Theorem \ref{res:main}.

Let $\rho>0$ and $V$, $\omega$
and $\Phi$ be as in the Theorem \ref{res:bourgain-jitomirskaya}.
As a consequence of Proposition \ref{res:zeromeasure}, the set
\[
A= \sigma^L(V,\omega)\setminus\bigcup_{\phi \not \in \Phi}
\sigma^L_{pp}(V,\omega,\phi),
\]
where $\sigma^L_{pp}(V,\omega,\phi)$ is the set of point eigenvalues 
of $L_{V,\omega,\phi}$ given by Theorem \ref{res:bourgain-jitomirskaya} satisfies that
$n_{V,\omega}^L(A)=0$. Indeed, according to Proposition \ref{res:zeromeasure} we only need to
show that $\mu_{\phi}(A)=0$ for all $\phi \not\in \Phi$, where
$\mu_{\phi}$ are the spectral measures of the
long-range operators $L_{V,\omega,\phi}$.
This is a consequence of the fact that the spectral
measures $\mu_{\phi}$, for $\phi \not\in \Phi$ are supported on
the set of point eigenvalues of the corresponding operator.

Therefore also $n^H_{V,\omega}(A)=0$ due to Proposition \ref{res:zeromeasure}. To 
prove Corollary \ref{res:bloch}  it only remains to show that also 
the Lebesgue measure of $A$ is zero. To do so, one can 
invoke Deift \& Simon \rmcite{deift-simon}. For almost periodic 
discrete Schr\"odinger operators they prove that for
Lebesgue almost every $a$ in
the set where the Lyapunov exponent is zero, one has the inequality
\begin{equation}\label{eq:inequality}
2 \pi\sin{\pi k^H_{V,\omega}(a)} \frac{dk^H_{V,\omega}}{da} \ge 1.
\end{equation}
Thus, under the additional assumption that the Lyapunov exponent
vanishes in the spectrum, the inequality (\ref{eq:inequality}) implies
that if $A$ is a subset of $\sigma^H(\omega,V)$ with
$n^H(A)=0$ then also the Lebesgue measure of $A$ is
zero.

As a consequence of Bourgain \& Jitormirskaya 
\rmcite{bourgain-jitomirskaya:absolutely,bourgain-jitomirskaya}, 
for any $a$ in $\sigma^H(V,\omega),$ (with
$|V|_{\rho}<\pet$) the Lyapunov exponent is zero.
 Therefore, the set $A$ has Lebesgue
measure zero and for the values of $a$ in its
complement in the spectrum,
\[
a \in \sigma^H_{V,\omega} \setminus A,
\]
which is a total measure subset of
$\sigma^H(V,\omega)$, the corresponding Harper-like equation
\begin{equation}\label{eq:harper-proof}
x_{n+1}+x_{n-1} + V(2\pi\omega n) x_{n} = a x_n, n \in \mZ
\end{equation}
has an analytic quasi-periodic Bloch wave, using the
argument of duality in the beginning of Section
\ref{sec:aubry}. Indeed, we saw that if $a$ is a
point eigenvalue of the operator $L_{V,\omega,\phi}$
whose eigenfunction decays exponentially then the
Harper-like equation (\ref{eq:harper-proof}) has an analytic quasi-periodic
Bloch wave with Floquet exponent $\phi$. This
completes the proof of Corollary \ref{res:bloch}.

\subsection{From Bloch waves to
reducibility}\label{sec:bloch2red}

In this section we will see how Corollary \ref{res:bloch} (which we proved in
the previous section) implies our main result, Theorem \ref{res:main}. By this
corollary we know that if $V$, $\omega$ and $\Phi$ are as in Theorem
\ref{res:bourgain-jitomirskaya} then, for almost all $a \in
\sigma^H(V,\omega)$, the equation (\ref{eq:harper-proof}) has
an analytic quasi-periodic Bloch wave with Floquet exponent $\fhi
\not\in \Phi$. Since we only want to prove a result for almost
every $a$, it is
sufficient to show that if $\fhi \not\in \Phi$ is such that
\begin{equation}\label{eq:nonresonant}
\fhi -\pi k \omega - \pi j \ne 0
\end{equation}
for all $k,j \in \mZ$ and (\ref{eq:harper-proof}) has an
analytic quasi-periodic Bloch wave with this Floquet exponent
$\fhi$, then the corresponding Schr\"odinger cocycle
$(A_{a,V},\omega)$ is reducible to constant coefficients.

\begin{remark}
If $\fhi/2\pi$ is resonant with respect to $\omega$,
\[
\fhi= \pi k + \pi j \omega,
\]
for some integers $k,j$, then one can also prove reducibility \rmcite{puig}. In
Section \ref{sec:cantor} we will consider the case of $\fhi=2\pi k$ which will
be used for the Cantor structure of the spectrum.
\end{remark}

The existence of a Bloch wave for Equation (\ref{eq:harper-proof}) 
implies that the Schr\"odinger cocycle has the following quasi-periodic 
solution 
\begin{equation}\label{eq:qpsolution}
\left(
\begin{array}{c}
{\tilde \psi}( 4\pi \omega  + \theta) \\
e^{-\ii\fhi}{\tilde \psi}( 2\pi \omega  + \theta)
\end{array}
\right)= e^{-\ii\fhi} 
\left( 
\begin{array}{cc}
a - V(\theta) & \; -1 \\
1             & \; \fm 0 
\end{array}\right)
\left(\begin{array}{c}
{\tilde \psi}( 2\pi \omega  + \theta) \\ 
{e^{-\ii\fhi}}{\tilde \psi}(  \theta)
\end{array}
\right)
\end{equation}
for all $\theta \in \mT$. Moreover, writing
\begin{equation}\label{eq:definiciov}
v(\theta)= \left({\tilde \psi}(\theta+2\pi\omega),{e^{-\ii\fhi}}{\tilde
\psi}(\theta)\right)^T
\end{equation}
and
\begin{equation}\label{eq:definicio-Y}
Y(\theta)= 
\left(\begin{array}{cc}
v_1(\theta) & \; {\bar v}_1(\theta) \\
v_2(\theta) & \; {\bar v}_2(\theta)
\end{array}\right),
\end{equation}
where the bar denotes complex conjugation, one always has the relation
\begin{equation}\label{eq:conjugationY}
A_{a,V}(\theta) Y(\theta)= Y(\theta+2\pi\omega) \Lambda(\fhi),
\end{equation}
where
\[
\Lambda(\fhi)=\left(
\begin{array}{cc}
e^{\ii\fhi} & 0 \\
0 & e^{-\ii\fhi}
\end{array}\right).
\]
Obviously, $Y$ will only define a conjugation between the cocycles
$(A_{a,V},\omega)$ and $(\Lambda(\fhi),\omega)$ if it is
nonsingular. Because of (\ref{eq:conjugationY}),
the determinant of $Y$ is constant as a function of $\theta$ and
it is  purely imaginary. In particular,  $v(\theta)$ and
${\bar v}(\theta)$ are linearly  independent for all $\theta$ if,
and only if, they are independent for some $\theta$. In the case
that $v$ and ${\bar v}$ are linearly independent, it is
not difficult to prove reducibility to constant coefficients of the
cocycle.

\begin{lemma}\label{res:reducibility_noress}
Let $A : \mathbb{T} \to SL(2,\mathbb{R})$ be a real analytic
map  and $\omega$ be nonresonant. Assume that there is an
 analytic map $v:\mathbb{T} \to \mathbb{R}^2\setminus\{0\}$, with $v$
and $\bar v$ linearly independent, such that
\[
v(\theta + 2\pi \omega) = e^{-\ii\fhi} A(\theta) v(\theta)
\]
holds for all $\theta \in \mathbb{T}$, where $\fhi \in [0,2\pi)$. 
Then the cocycle $(A,\omega)$ is reducible to constant
coefficients by means of a real analytic transformation.
Moreover, the Floquet matrix can be chosen to be
of the form
\begin{equation}\label{eq:floquetformirr}
B= \left(
\begin{array}{cc}
\fm\cos{\fhi} & \; \sin{\fhi} \\
-\sin{\fhi} & \;   \cos{\fhi}
\end{array}
\right).
\end{equation}
\end{lemma}

{\noindent{\bf{Proof:}}} Let $Z^1(\theta)= Y(\theta)$ as in
(\ref{eq:definicio-Y}), $B^1=\Lambda(\fhi)$ and
\[
d(\theta)= v_1(\theta) {\bar v}_2(\theta) 
- {\bar v}_1(\theta)v_2(\theta).
\]
be the determinant of $Z^1$. Therefore $Z^1$ defines a conjugation between
$(A_{a,V},\omega)$ and $(B^1,\omega)$ because $v$ and $\bar v$ are
linearly independent, $Z^1$ is real analytic and, 
for every $\theta \in \mT$, $Z^1(\theta)$ is nonsingular.

Moreover, from the conjugacy (\ref{eq:conjugationY}) and the
nonresonance of $\omega$,  $d(\theta)$ is constant as a function of
$\theta$. By the linearity of our system, we choose this constant
value to be $-\ii/2$ (recall that, due to the form of $Z_1,$ its
determinant must be purely imaginary).

To obtain the real rotation consider the composition
\[
Z(\theta)=  Z^1(\theta) Z^2
\]
where $Z^2$ is the constant matrix
\[
Z^2=
\left(
\begin{array}{cc}
1 & -\ii \\
1 & \fm\ii
\end{array}\right),
\]
Then $Z$ satisfies the desired conjugation
\[
A(\theta) Z(\theta)= Z(\theta+2\pi\omega) B
\]
being $B$ the rotation of angle $\fhi$ given by (\ref{eq:floquetformirr}). Thanks to the
construction $Z$ is real and with determinant one.
\hfill\QED

To complete the proof of Theorem \ref{res:main}, it only remains to rule 
out the possibility that $\fhi$
satisfies (\ref{eq:nonresonant}) and 
 $v$ and ${\bar v}$ are linearly dependent at the same
time. Recall that these two vectors are linearly independent for all
$\theta$ if, and only if, they are linearly independent for some $\theta.$
Note that both $v(\theta)$ and ${\bar v}(\theta)$ are different
from zero for all $\theta \in \mT$ by construction. Assume that
$v(\theta)$ and ${\bar v}(\theta)$ were linearly dependent for all
$\theta$. Since these vectors depend analytically on $\theta$,
there would exist an analytic $h:\mT \to \mR$ and an integer $k \in
\mZ$ such that
\[
{\bar v}(t)= e^{\ii\left(h(t)+kt\right)} v(t)
\]
for all $t \in \mR$. Using that $v$ and ${\bar v}$ are
quasi-periodic solutions of $(A,\omega)$, this would imply that
\[
e^{\ii\left(h(t)+kt\right)} e^{\ii\fhi}= 
e^{\ii\left(h(t+2\pi\omega)+kt+2\pi k \omega\right)} e^{-\ii\fhi}.
\]
Therefore, $h$ must satisfy the
following small divisors equation\index{small divisors equation}
\[
h(\theta+2\pi\omega)-h(\theta)= 2\fhi - 2\pi k\omega -2\pi j
\]
for all $\theta \in \mT,$ where $j$ is some fixed integer. Clearly, 
such analytic $h$ cannot exist unless
\[
\fhi= \pi \left(j + k \omega \right),
\]
which is a contradiction with the nonresonance condition
(\ref{eq:nonresonant}). This ends the proof of Theorem
\ref{res:main}. \hfill\QED

\section{Applications}\label{sec:applications}

In this section we will prove several consequences of the main theorem
which are summarized in theorems \ref{res:cantor} and
\ref{res:nored}. In Section \ref{sec:gap_edges} we will present the setting
of this section. In  \ref{sec:moser-poschel} we adapt
Moser-P\"oschel perturbation arguments to the discrete
case. This is applied in Section \ref{sec:cantor} to
the proof of nonperturbative genericity of Cantor spectrum.
Finally, in Section \ref{sec:nonred} we
prove that, in our situation, Cantor
spectrum implies nonreducibility for a $G_\delta$-set of
energies.

\subsection{Reducibility at gap edges}\label{sec:gap_edges}

In previous sections, we discussed the reducibility of a quasi-periodic
Schr\"odinger cocycle $(A_{a,V},\omega)$ when $a$ is a point eigenvalue of
the dual operator $L_{V,\omega,\phi}$ and $\phi$ satisfies a nonresonance
condition of the form (\ref{eq:nonresonant}), which was enough to prove
the main result. The ``resonant'' values of $\phi$:
\begin{equation}\label{eq:resonantphis}
\phi= \pi j + \pi\omega k, \qquad j,k \in \mZ
\end{equation}
are particularly important for the description of the spectrum of
these operators because the corresponding point eigenvalues lie at
endpoints of spectral gaps. Let us prove the reducibility at these
endpoints. What follows mimics the proof of the ``Ten Martini Problem''
given in \rmcite{puig}.

Bourgain \& Jitomirskaya \rmcite{bourgain-jitomirskaya:absolutely} also
prove that, provided $|V(\theta)| <\pet$ and $\omega$ is 
Diophantine $L_{V,\omega,\phi}$ has pure-point spectrum with 
exponentially localized eigenfunctions if $\phi$ is of the form
(\ref{eq:resonantphis}), see Remark 8.2 after Theorem 7 in
\rmcite{bourgain-jitomirskaya:absolutely}. Taking into
account the symmetries of the operators, the study reduces to the four cases
$\phi=0,\pi,\pi\omega, \pi\omega +\pi$. For the sake of simplicity we
consider here $\phi=0$. As a direct consequence of Aubry duality and
Bourgain-Jitomirskaya result for the dual operators $L_{V,\omega,0}$,
the set of pure-point eigenvalues $\sigma_{pp}^L(V,\omega,0)$ is a dense
subset of $\sigma^H(V,\omega)$ and any energy $a\in \sigma_{pp}^L(V,\omega,0)$ has
a quasi-periodic Bloch wave with $0$ as a Floquet exponent for the dual eigenvalue equation. Let $a$ be one
of these eigenvalues. The condition above means that there is an analytic
map $\tilde{\psi}:\mT \to \mC$ such that
\[
x=(x_n)_{n \in \mZ} = \left(\tilde{\psi}\left(2\pi\omega n +
\theta\right)\right)_{n \in \mZ}
\]
is a nonzero solution of $H_{V,\omega,\theta} x= a x$. Clearly, due to the
symmetry of the eigenvalue equation for $L_{V,\omega,0}$, the function
$\tilde{\psi}$ can be chosen real analytic. In terms of the cocycle we
have that the relation 
\[
\left(
\begin{array}{c}
{\tilde \psi}( 4\pi \omega  + \theta) \\
{\tilde \psi}( 2\pi \omega  + \theta)
\end{array}
\right)=\left( 
\begin{array}{cc}
a - V(\theta) & \; -1 \\
1                                        & \; \phantom{-} 0 
\end{array}\right)
\left(\begin{array}{c}
{\tilde \psi}( 2\pi \omega  + \theta) \\
{\tilde \psi}(  \theta)
\end{array}
\right)
\]
holds for all $\theta\in\mT$. Instead of Lemma \ref{res:reducibility_noress}
we now have the following (see \rmcite{puig} for the proof,
which is a simple triangularization and averaging argument).

\begin{lemma}\label{res:reducibility_ress}
Let $A \in C^a_\delta(\mathbb{T},SL(2,\mathbb{R}))$ be a
real analytic map and $\omega$ be Diophantine. 
Assume that there is a nonzero real analytic map $v \in
C^a_\delta(\mT,\mR^2)$ such that the relation
\[
v(\theta + 2\pi \omega) = A(\theta) v(\theta)
\]
holds for all $\theta \in \mathbb{T}.$ Then the quasi-periodic 
cocycle $(A,\omega)$ is reducible to constant coefficients 
by means of a quasi-periodic transformation which
is analytic in $|\Im \theta| < \delta$.
 Moreover the Floquet matrix can be chosen to be of the form
\begin{equation}\label{eq:floquetform}
B= \left(
\begin{array}{cc}
1 & \; c \\
0 & \; 1 
\end{array}
\right)
\end{equation}
for some $c \in \mR$.
\end{lemma}

In the Almost Mathieu case Ince's argument
\rmcite{ince, puig} shows
that $c \ne 0$. Otherwise the dual model (which is also a
Schr\"odinger operator) would have a point eigenvalue with two
linearly independent eigenvectors in $l^2(\mZ)$, and this is
a contradiction with the limit-point character of
Schr\"odinger operators (or just the preservation of the
Wronskian in this discrete case). The fact that $c\ne 0$ is
important for the description of the spectrum, because if a
Schr\"odinger cocycle is reducible to a Floquet
matrix of the form (\ref{eq:floquetform}) with $c\ne 0$
then  the corresponding energy is at the endpoint
of an open gap in the spectrum, as it will
be seen in the next section. 

For general potentials $V$, however, we cannot use Ince's argument and it
may happen that some of these are collapsed. In fact, there
are examples of quasi-periodic Schr\"odinger operators (with $V$ small,
$\omega$ Diophantine) for which some $c$ are zero
\rmcite{broer-puig-simo} or even do not display
Cantor spectrum (see De Concini 
\& Johnson \rmcite{deconcini-johnson}). 

Nevertheless, even if $c$ can be zero, Moser \& P\"oschel
\rmcite{moser-poschel} showed that, in this reducible
setting, a closed gap can be opened by means of an
arbitrarily small and generic real analytic perturbation 
of the potential. In the next
section we give an adaption of their proof to the discrete
case together with some extra properties which will be needed later.

\subsection{Moser-P\"oschel perturbation
argument}\label{sec:moser-poschel}

In this section we prove the following adaption of Moser-P\"oschel
argument to the discrete case, which deals with cocycles
which are perturbations of constant matrices of the form
(\ref{eq:floquetform}).

\begin{proposition}\label{res:perturbation}
Let $V$ be real analytic, $\omega$ Diophantine and $(A_{a,V},\omega)$, for
some $a \in \mR$ be a quasi-periodic Schr\"odinger cocycle. Assume that $(A_{a,V},\omega)$ 
is analytically reducible to the constant coefficients cocycle $(B,\omega)$
with
\[
B=\left(
\begin{array}{cc}
1 & c \\
0 & 1 
\end{array}\right)
\]
for some $c \in \mR$. Let $W:\mT \to \mR$ be
real analytic and $\alpha$ real. If $Z:\mT \to SL(2,\mR)$ is the real
analytic reducing matrix and the conditions 
\begin{equation}\label{eq:conditioncne0}
c\ne 0 \quad\text{ and }\quad [Wz_{11}^2] \ne 0
\end{equation}
or 
\begin{equation}\label{eq:conditionce0}
c= 0 \quad\text{ and }\quad -[Wz_{11}z_{12}]^2+[W z_{12}^2][Wz_{11}^2] > 0
\end{equation}
are satisfied, then the quasi-periodic cocycle $(A_{a,V+\alpha W},\omega)$
has an exponential dichotomy provided $|\alpha|>0$ is small enough and 
\begin{equation}\label{eq:conditionalpha}
c\alpha[Wz_{11}^2] < 0 \text{ if  (\ref{eq:conditioncne0})
holds.}
\end{equation}
Moreover in the case $c \ne 0$, the Lyapunov exponent of
$(A_{a,V+\alpha W},\omega),$ $\gamma(a,V+\alpha W,\omega),$ and its
rotation number, $\rot_f(a,V+\alpha W,\omega),$ satisfy   
\begin{equation}\label{eq:asymptoticsW}
\lim_{\tiny{\begin{array}{c}{\alpha \to 0,}\\{c\alpha[Wz_{11}^2] >
0}\end{array}}} \frac{|\gamma(\alpha)|}{|\alpha|^{1/2}}=
c [Wz_{11}^2]
=
\lim_{\tiny{\begin{array}{c}{\alpha \to
0,}\\{c\alpha[Wz_{11}^2] < 0}\end{array}}}
\frac{|\rot_f(\alpha)-\rot_f(0)|}{|\alpha|^{1/2}}
\end{equation}
\end{proposition}

\begin{remark}
This type of perturbation arguments have been used in a
variety of contexts, c.f.  Moser \& P\"oschel
\rmcite{moser,moser-poschel},  Johnson \rmcite{johnson:cantor},
 N\'u\~nez \rmcite{nunez}, Broer, Puig \& Sim\'o
\rmcite{broer-puig-simo} and Puig \& Sim\'o \rmcite{puig-simo}.
\end{remark}

\noindent{\bf{Proof:}} Since $Z$ is the reducing transformation of 
$(A_{a,V},\omega)$ to $(B,\omega)$, it also renders the perturbed cocycle
$(A_{a,V+\alpha W},\omega)$ to $(B+ \alpha W P,\omega)$
where
\[
P(\theta)=
\left(
\begin{array}{cc}
 z_{11} z_{12} - c z_{11}^2        &    
 -c z_{11} z_{12} + z_{12}^2 \\
- z_{11}^2  & \;\;\;  - z_{11} z_{12}
\end{array}
\right).
\]

After one step of averaging this cocycle can be analytically conjugated to 
\[
\left(B + \alpha [W P] + \alpha^2 R_2,\omega \right)
\]
where $[\cdot]$ denotes the average of a quasi-periodic
function and $R_2$ depends analytically on $\alpha$ and
$\theta$ in some open neighbourhoods of $0$ and $\mT$.
Moreover, a computation shows that
\[
B + \alpha [W P] + \alpha^2 R_2 = \exp\left(\tilde{B}_0 + \alpha \tilde{B}_1
+ \alpha^2 \tilde{R_2}\right),
\]
being
\[
\tilde{B}_0=\left(
\begin{array}{cc}
0 &  c \\
0 & 0  
\end{array}\right), \quad
\tilde{B}_1 = \left(
\begin{array}{cc}
 [W z_{11} z_{12}] - \frac{c}{2} [W z_{11}^2] &  -c [W z_{11} z_{12}] + [W z_{12}^2] \\
-[W z_{11}^2]  &   - [W z_{11} z_{12}] + \frac{c}{2} [W z_{11}^2]
\end{array}\right)
\]
and $\tilde{R_2} \in sl(2,\mR)$ depending analytically on $\alpha$ and
$\theta$. Let 
\[
D= \left(\begin{array}{cc}
d_1 & \fm d_2 \\
d_3 & -d_1
\end{array}\right)
=\tilde{B}_0 + \alpha\tilde{B}_1,
\]
whose determinant is $d= -d_1^2 -d_2 d_3$. Now let us distinguish
between the cases $c\ne 0$  and $c=0$.

If (\ref{eq:conditioncne0}) holds then the expression for
the determinant becomes 
\[
d= c \alpha [W z_{11}^2] +O(\alpha^2)
\]
so that it is negative if, in addition
(\ref{eq:conditionalpha}) holds. In this case, the matrix
\[
Q= \left(\begin{array}{cc}
d_2 & d_2 \\
-d_1 + \sqrt{-d} & - d_1 - \sqrt{-d}
\end{array}\right),
\]
which is well-defined, has  determinant
\[
2c \sqrt{-c \alpha [W z_{11}^2]} + O(\alpha).
\]
and satisfies $DQ= Q \Delta,$ where
\[
\Delta= \left(\begin{array}{cc}
\sqrt{-d} & 0 \\
0 & -\sqrt{-d}.
\end{array}\right).
\]

Therefore the change of variables defined by $Q$ transforms the
cocycle 
\begin{equation}\label{eq:pompeu1}
\left(\exp\left(\tilde{B}_0 + \alpha \tilde{B}_1 + \alpha^2
\tilde{R_2}\right),\omega \right)
\end{equation}
into 
\begin{equation}\label{eq:pompeu2}
\left(\exp\left(\left(\Delta+\tilde{S}_{2}\right)\right),\omega \right)
\end{equation}
where
\[
\tilde{S}_2(\alpha,\theta) = \alpha^2 Q^{-1} R_2(\alpha,\theta) Q 
\]
which is $O(|\alpha|^{3/2})$ uniformly in $\theta$. Note that
\[
\Delta+\tilde{S}_{2} =
\sqrt{-c \alpha [W z_{11}^2]} \left(\left(\begin{array}{cc} 1 & \fm 0 \\ 0 & -1
\end{array}\right) +O(|\alpha|) \right).
\]
so that if (\ref{eq:conditionalpha}) holds and $|\alpha|>0$ is small
enough the cocycle $(A_{a,V+\alpha W},\omega)$ has an exponential
dichotomy and the Lyapunov exponent satisfies
\[
\frac{1}{\sqrt{|\alpha|}} \left|\gamma\left(A_{a,V+\alpha
W},\omega\right)\right| \to \sqrt{|c [W z_{11}^2]|},
\]
see Coppel \rmcite{coppel}. To obtain the asymptotics of the
rotation number, we can consider the transformation $Q$
defined for $c \alpha [W z_{11}^2]>0$. This, although
complex, is a well-defined conjugation between $D$ and $\Delta$,
which is now a complex rotation of angle $\sqrt{|d|}.$
Therefore (\ref{eq:pompeu1}) is conjugated to
(\ref{eq:pompeu2}), a perturbation of a complex rotation.
Using the definition of the fibered rotation number given in
Section \ref{sec:eliasson} the result follows.

Let us now consider the situation when
(\ref{eq:conditionce0}) holds. In 
this case, the matrix $D$ becomes
\[
D= \alpha \left(
\begin{array}{cc}
 [W z_{11} z_{12}]   &   [W z_{12}^2] \\
-[W z_{11}^2]  &   - [W z_{11} z_{12}] 
\end{array}\right)=: \alpha \tilde{D}
\]
Condition (\ref{eq:conditionce0}) is equivalent to the hyperbolicity
of $\tilde{D}$ whose determinant is
\[
\tilde{d} = - [W z_{11} z_{12}]^2 + [W z_{12}^2] [W z_{11}^2].
\]
Therefore there is a change of variables $Q$, independent of $\alpha$ and
$\theta,$ which renders it to a diagonal form $\tilde{\Delta}$ with
$\sqrt{-d}$ and $-\sqrt{-d}$ as dia\-gonal entries. This conjugation
transforms the cocycle 
\[
\left(\exp\left(\alpha \tilde{D} + \alpha^2
\tilde{R_2}\right),\omega \right)
\]
into 
\[
\left(\exp\left(\left(\alpha\tilde{\Delta}+\alpha^2\tilde{S}_{2}\right)\right),\omega \right)
\]
where
\[
\tilde{S}_2(\alpha,\theta) =  Q^{-1} \tilde{R}_2(\alpha,\theta) Q.
\]
Since, 
\[
\alpha\Delta+\alpha^2\tilde{S}_{2} =
\alpha \sqrt{-\tilde{d}} \left(\left(\begin{array}{cc} 1 & \fm 0 \\ 0 & -1
\end{array}\right) + \frac{\alpha}{\sqrt{-\tilde{d}}}\tilde{S}_2 \right).
\]
the cocycle $(A_{a,V+\alpha W},\omega)$ 
has an exponential dichotomy when $\alpha \ne 0$ is small
enough  (see, again Coppel \rmcite{coppel}). \hfill\QED

The perturbation argument in the previous proposition can applied to the
reducible cocycles at endpoints of gaps as we do next.

\begin{corollary}\label{res:ncg}
Let $V$, $a$ and $\omega$ as in Proposition \ref{res:perturbation} and
assume that $c \ne 0$. Then $a$ is at the endpoint of a noncollapsed
spectral gap $I$ of $\sigma(V,\omega)$ (the right one if
$c>0$ and the left one if $c<0$). Moreover, the limits 
\begin{equation}\label{eq:asymptoticsncg}
\lim_{\tiny{\begin{array}{c}{\alpha \to 0,}\\{a+\alpha \in I}\end{array}}} 
\frac{\gamma(a+\alpha,V,\omega)}{\sqrt{|\alpha|}}=
\lim_{\tiny{\begin{array}{c}{\alpha \to 0,}\\{a+\alpha \in I}\end{array}}} 
\frac{\left|\rot_f(a+\alpha,V,\omega)-\rot_f(a,V,\omega)\right|}{\sqrt{|\alpha|}}
\end{equation}
exist and are different from zero.
\end{corollary}

\noindent{\bf{Proof:}} Take $W=1$ in Proposition \ref{res:perturbation}.
Then, the cocycle $(A_{a+\alpha,V},\omega)$ has an exponential dichotomy
if $c\alpha<0$ and $|\alpha|$ is small enough. This means that there is an
open spectral gap besides $a$ (to the left if $c>0$ and to the right
otherwise). Moreover the asymptotics of formula
(\ref{eq:asymptoticsW}) imply (\ref{eq:asymptoticsncg}). \hfill\QED

Finally we consider the variation of the rotation number in the case $c=0$ in a more general
setting which will be needed in the next section. 

\begin{proposition}\label{res:asymptotics_nogap}
Let $V$ be continuous and $\omega$ nonresonant. Assume that the
Schr\"odinger cocycle $(A_{a_0,V},\omega)$ is reducible to
the cocycle $(B,\omega)$ with $B \in SO(2,\mR)$ a constant
matrix. Then the map
\[
a \in \mR \mapsto \rot_f(a,V,\omega)
\]
is differentiable at $a_0.$
\end{proposition}

\noindent{\bf{Proof:}} Let $\rho$ be the angle of the
rotation,
\[
B=\left(
\begin{array}{cc}
\cos{\rho} & \sin{\rho} \\
-\sin{\rho} & \cos{\rho}
\end{array}
\right).
\]
The cocycle $(A_{a,V},\omega)$ is conjugated to
$(B+\alpha R,\omega)$ where 
\[
R(\theta)= Z(\theta+2\pi\omega)^{-1} 
\left(\begin{array}{cc}
1 & 0 \\
0 & 0
\end{array}\right)
Z(\theta)
\]
and $\alpha= a-a_0$. The cocycle $(B+\alpha R,\omega)$
induces a lift  $\tilde{F}$ from $\mR \times \mT$ to itself
of the form
\[
{\tilde F}(t,\theta)=\left(t+\rho + \alpha f(t,\theta,\alpha),
\theta+2\pi\omega\right),
\]
where $f$ is continuous and $2\pi$-periodic in both $t$ and
$\theta$. Therefore,
\begin{multline}
\rot_f(B+\alpha R,\omega)-\rot_f(B,\omega)= \\
\lim_{N \to
\infty} \frac{1}{2\pi N} \sum_{n=0}^{N-1} \left(
\alpha f\left(\tilde{F}^n(t,\theta)\right)\right)= O(\alpha)
\nonumber
\end{multline}
as we wanted to show. \hfill\QED

\begin{remarks}
\mbox{}
\begin{enumerate}
\item A computation shows that the derivative of the rotation number above is nonzero. In particular, when a Schr\"odinger cocycle is reducible to the identity, the corresponding energy lies at the endpoint of a collapsed gap.       
\item Similar results have been obtained when $Z,Z^{-1}:\mT \to SL(2,\mR)$ are
square integrable and $B:\mT \to SO(2,\mR)$ is measurable, compare with Moser \rmcite{moser} and Deift \& Simon \rmcite{deift-simon}. 
\end{enumerate}
\end{remarks}

\subsection{Genericity of Cantor spectrum}\label{sec:cantor}

In the previous section we have seen that if a Schr\"odinger
cocycle is reducible to a matrix with trace $2$ then the
corresponding energy is at the endpoint of a spectral gap which is collapsed if the Floquet matrix is the identity. The next consequence  of Proposition
\ref{res:perturbation} is that when the Floquet matrix is the identity (a similar statement holds for minus the identity) one can ``open up''  the collapsed gap by means of a generic perturbation.

\begin{corollary}\label{res:cg}
Let $V$, $a$, $\omega$ and $Z$ be as in Proposition \ref{res:perturbation} and
assume that $c = 0$. If $W$ is a generic real analytic potential
then for $|\beta|\ne 0$ small enough the spectrum $\sigma(V+\beta
W,\omega)$ has an open spectral gap with  \IDS $k(a,V,\omega)$.
\end{corollary}

\noindent{\bf{Proof:}}  In Proposition \ref{res:perturbation} we proved
that for a perturbation $\tilde{W}$ satisfying (\ref{eq:conditionce0})
\[
-[\tilde{W}z_{11}z_{12}]^2+[\tilde{W} z_{12}^2][\tilde{W}z_{11}^2] < 0
\]
the cocycle $(A_{a,V+\beta \tilde{W}},\omega)$ has an exponential
dichotomy if $|\beta|>0$ small enough. This means that $a$ lies in a spectral gap of
$\sigma(V+\beta \tilde{W},\omega)$  which, by continuity must satisfy that
\[
k(a,V+\beta \tilde{W},\omega)= k(a,V,\omega)
\]
for $|\beta|$ small enough. Let us now show that if $W$ is a
generic potential, then for every $|\beta|\ne 0$ small enough there is a value of
$\alpha$ such that $a+\alpha$ lies in a spectral gap of $\sigma(V+\beta
W,\omega)$ with
\[
k(a+\alpha,V+\beta W,\omega)= k(a,V,\omega).
\]
Note that the condition (\ref{eq:conditionce0}) can be rewritten as
\[
[{\tilde W} y_1]^2 - [\tilde{W} y_2]^2 - [\tilde{W} y_3]^2 <0
\]
where
\[
y_1= \frac{1}{2}\left(z_{11}^2 + z_{12}^2\right),\quad
y_2= \frac{1}{2}\left(z_{11}^2 - z_{12}^2\right),\quad
y_3= z_{11}z_{12}.
\]
Let $\alpha$ be such that
\[
[(\alpha +W) y_1]= \alpha[y_1]+[W y_1] = 0,
\]
(this determines $\alpha$ since $[y_1]\ne 0$). 
Then the shifted perturbation $\alpha +W$ 
satisfies condition (\ref{eq:conditionce0}) unless
\[
-[W y_1] [y_2] +  [W y_2] [y_1]= 0 
\quad\text{and}\quad
-[W y_1] [y_3] +  [W y_3] [y_1]= 0,
\]
which is clearly a generic condition. Then, if
$|\beta|>0$ is small enough, the spectrum $\sigma(V+\beta W,\omega)$
 has an open gap with
\[
k(a+\alpha\beta,V+\beta W,\omega)= k(a,V,\omega)
\]
as we wanted to show. \hfill\QED

\begin{remark}
As Moser \& P\"oschel show, when $c=0$ it is always possible to choose the
reducing transformation such that $[z_{11}^2]=[z_{12}^2]=1$ and
$[z_{11}z_{12}]=0$ so that $[y_1]= 1$, $[y_2]=0$, $[y_3]=0$ 
and the generic $W$ must satisfy
\[
[W (z_{11}+z_{12})^2 ] \ne 0  \qquad\text{or}\qquad [W
(z_{11}-z_{12})^2 ] \ne 0.
\]
\end{remark}

Let us now summarize the situation. Using the two past
sections we have seen that if $V$ is a real analytic potential on
$C^a_\rho(\mT,\mR)$, with $|V|_{\rho}<\pet$ and $\omega$ is Diophantine,
there is a countable dense subset of energies in the spectrum where the 
system is reducible to a Floquet matrix with trace $2$. 
These lie at endpoints of gaps. Although
these can be collapsed, a generic and arbitrarily small perturbation opens
them as Corollary \ref{res:cg} says. Since there is a countable
number of gaps Theorem \ref{res:cantor} follows.

\subsection{Cantor spectrum implies
nonreducibility}\label{sec:nonred}

In \rmcite{eliasson:floquet} it was seen that, for a generic
real analytic Schr\"odinger cocycle (with Diophantine frequencies) besides
the almost everywhere reducibility there was a set of zero measure of
energies for which the cocycle was not reducible to constant coefficients.
The proof relies on the {\sc{kam}} procedure developed there, but the
Cantor structure of the spectrum is seen to play a key role. In this
section we prove irreducibility for a $G_\delta$-set of energies assuming
only Cantor structure of the spectrum and Theorem
\ref{res:main}. This argument is reminiscent of some
techniques in circle maps, see Arnol'd
\rmcite{arnold:circle}. We state here a slightly more general version than
that of Theorem \ref{res:nored}. More applications will be
given elsewhere.

\begin{theorem}\label{res:nonreducible}
Let $\rho>0$. There is a constant $\pet>0$ such that if $V \in
C^a_\rho(\mT,\mR)$ is real analytic with $|V|_{\rho}<\pet$, $\omega$ is 
Diophantine and $I$ is an open interval such that 
\[
K= \sigma(V,\omega) \cap {\bar I}
\]
is a nonvoid Cantor set, then there is a $G_\delta$-dense set
of energies in $K$ for which the corresponding
Schr\"odinger cocycle is not reducible to constant
coefficients by means of a continuous transformation.
\end{theorem}

\noindent{\bf{Proof:}} Consider, for any 
$a_1,a_2 \in K$ with $a_1 \ne a_2$,
\[
\delta(a_1,a_2)= \left|
\frac{k(a_1,V,\omega)-k(a_2,V,\omega)}{a_1-a_2}\right|.
\]
Now, for any $a \in \mR$ we can define
\[
m(a)= \sup_{\lambda \ne a,\lambda \in K} \delta(a,\lambda)
\]
which is either a positive real number or $+\infty.$ 

If $a \in \sigma(V,\omega)$ is reducible to constant
coefficients then we have two situations. Either the Floquet
matrix $B$ has trace $\pm 2$, in which case $m(a)=\infty$ (see Corollary
\ref{res:ncg}) or $B \in SO(2,\mR)$ and then $m(a)<\infty$ (see Proposition
\ref{res:asymptotics_nogap}). Due to the fact that
$|V|_\rho< \pet,$ 
 $\omega$ is Diophantine and the Cantor structure of the spectrum there is
a dense set of endpoints of gaps, $\cG_K$, where the system is
reducible to constant coefficients because of Eliasson's
Theorem \ref{res:eliasson}.

We will show that the set where $m(a)=\infty$ is a
$G_{\delta}$-dense subset of $K$. Excluding the endpoints of
gaps where there is reducibility to a Floquet matrix with
trace $\pm 2$ (which are at most countable) we will still
have a $G_{\delta}$-dense subset of energies in $K$ whose
corresponding cocycle cannot be reducible to constant
coefficients.

Let, for any $n \in \mN \cup \{0\}$ and $a_0 \in K$,
\[
U(a_0,n)= \left\{a \in K; \delta(a,a_0) > n \right\}
\]
and 
\[
U(n)= \bigcup_{a_0 \in K} U(a_0,n).
\]
The sets $U(n)$ are open in $K$ because of the continuity of the
rotation number. Moreover they are dense in $K$ because they
contain $\cG_K$, which is dense in $K$. Therefore
\[
U(\infty) \bigcap_{n>0} U(n) = \left\{ a \in K; m(a) = \infty \right\},
\]
is a $G_\delta$-dense subset of $K$. If we exclude the endpoints of open
gaps  the remaining energies, which still form a
$G_\delta$-dense subset of $K$, cannot be reducible to constant
coefficients by means of a continuous transformation. This
proves \ref{res:nonreducible} and also \ref{res:nored}.
\hfill\QED

\section*{Acknowledgements}

The author wishes to thank H. Eliasson, S. Jitomirskaya, C.
Sim\'o and B. Simon for 
many ideas which have been decisive in this paper. The paper arose in a stay at
Caltech Mathematics Department and completed mostly in the \emph{Departament de
Matem\`atica Aplicada i An\`alisi} at the \emph{Universitat
de Barcelona} as part of his PhD thesis. I am grateful to both institutions. This work
has been supported by grants DGICYT BFM2003-09504-C02-01
(Spain) and  CIRIT 2001 SGR-70 (Catalonia).


\def\cprime{$'$}

\end{document}